\documentclass{amsart}
\usepackage{amsmath,amssymb,amscd,amsthm}

\newcommand{\C}[1]{\mathcal{#1}}
\newcommand{\G}[1]{\mathfrak{#1}}
\newcommand{\B}[1]{\mathbb{#1}}
\newcommand{\xra}[1]{\xrightarrow{\ #1\ }}
\newcommand{\fotimes}[1]{\ \underset{#1}{\Box}\ }
\newcommand{\eotimes}[1]{\underset{#1}{\otimes}}

\newtheorem{thm}{Theorem}[section]
\newtheorem{cor}[thm]{Corollary}
\newtheorem{lem}[thm]{Lemma}
\newtheorem{conj}[thm]{Conjecture}
\newtheorem{prop}[thm]{Proposition}

\theoremstyle{definition}
\newtheorem{defn}[thm]{Definition}
\newtheorem{rem}[thm]{Remark}
\newtheorem{exm}[thm]{Example}

\numberwithin{equation}{section}

\title{Excision in Hopf Cyclic Homology}
\author{Atabey Kaygun}
\email{{\tt akaygun@uwo.ca}}
\author{Masoud Khalkhali}
\email{{\tt masoud@uwo.ca}}
\address{Department of Mathematics\\
The University of Western Ontario\\
London, Ontario N6A 5B7\\
Canada}

\begin{document}
\maketitle

\section{Introduction}

Cyclic cohomology theory for Hopf algebras was introduced by Connes
and Moscovici in their ground breaking paper on transverse index
theory \cite{ConnesMoscovici:HopfCyclicCohomology} (cf. also
\cite{ConnesMoscovici:HopfCyclicCohomologyIa,ConnesMoscovici:HopfCyclicCohomologyII},
and \cite{ConnesMoscovici:HopfRelativeCyclic} for a recent survey and
new results).  This theory can be regarded as the right noncommutative
analogue of Lie algebra homology because of existence of a
characteristic map for actions of Hopf algebras and since it reduces
to Lie algebra homology for Hopf algebras associated to Lie algebras.
This should be compared with the role played by ordinary cyclic
cohomology of noncommutative algebras and its reduction to de Rham
homology for smooth commutative algebras
\cite{Connes:NonCommutativeGeometry}. The algebraic underpinnings of
Hopf cyclic theory however turned out to be more complicated. It is
now understood as a special case of a general theory developed in the
series of papers \cite{Khalkhali:DualCyclicHomology,
Khalkhali:InvariantCyclicHomology, Khalkhali:HopfCyclicHomology,
Khalkhali:SaYDModules} (with the help of
\cite{KhalkhaliAkbarpour:EquivariantCyclicHomology}) for (co)algebras
endowed with a (co)action of a Hopf algebra and where, among other
things, the right notion of coefficients was introduced (cf. also
\cite{Khalkhali:IntroHopfCyclicCohomology} for a recent survey). This
theory was later extended by the first author
\cite{Kaygun:BialgebraCyclicK} to bialgebras and to a much larger
class of coefficients.  This extension is rather surprising since
existence of a bijective antipode was essential in constructing the
cyclic and cocyclic complexes developed for Hopf cyclic homology of
both variants.  By bialgebra or Hopf cyclic homology we now mean the
theory expounded in these papers mentioned above.

The main results of the present paper are two excision theorems
(Theorem~\ref{Excision1} and Theorem~\ref{Excision2} below) for Hopf
cyclic homology of module coalgebras and comodule algebras
respectively.  A basic property of cyclic (co)homology for algebras
and coalgebras is excision. If by cyclic theory we mean the positively
graded cyclic theory then this question is completely settled by
Wodzicki's excision theorem \cite{Wodzicki:Excision}. Originally
formulated for algebras, this result states that a not necessarily
unital $k$-algebra $I$ satisfies the excision property in cyclic
homology for all algebra extensions $0\rightarrow I \rightarrow A
\rightarrow A/I \rightarrow 0$ if and only if $I$ is $H$-unital
(cf. \cite{FarinatiSolotar:ExcisionInCyclicHomologyOfCoalgebras} for
the corresponding coalgebra result). It is this result that we
generalize in this paper by finding sufficient conditions for excision
to hold in Hopf cyclic cohomology. It is known after Cuntz and Quillen
\cite{QuillenCuntz:Excision} that for periodic cyclic cohomology
excision always holds irrespective of $I$ being $H$-unital or
not. Wodzicki's theorem is used in the proof of this result and we
believe that the corresponding statements hold for periodic Hopf
cyclic homology as well and we state two conjectures to this effect
(Conjecture ~\ref{Conjecture1} and Conjecture ~\ref{Conjecture2}
below).

Throughout this paper we will work over a field $k$ of an arbitrary
characteristic.  One can prove all results we prove in this paper by
requiring all $k$--modules we use to be $k$--projective or more
generally $k$--flat, if one needs $k$ to be just a commutative algebra
with a unit.  Unadorned tensor products $\otimes$ are over $k$ and all
other tensor products will be clearly indicated.  All algebras
(resp. coalgebras) will be associative (resp.  coassociative) but not
necessarily unital (resp. counital).  The comultiplication and counit
of a coalgebra $C$ will be denoted by $\Delta$ and $\varepsilon$,
respectively.  The unit of an algebra $A$ will be denoted by
$\B{I}$. $B$ will always be a unital and counital, associative and
coassociative bialgebra, or a Hopf algebra.  All modules over $B$ are
assumed to be left modules unless it is stated otherwise. For a
$B$--module $X$ we will distinguish an element $b\in B$ and the
$k$--module endomorphism $L_b$ of $X$ defined by the left action of
$b$ on $X$.  If $X$ is simultaneously a $B$--module and a
$B$--comodule with no compatibility between these structures assumed,
we will call $X$ a $B$--module/comodule.  We will use Sweedler's
notation denoting comultiplication structure on coalgebras by
$\Delta(c) = c_{(1)}\otimes c_{(2)}$.  Similarly if $X$ is a left
(resp. right) comodule over a coalgebra $C$ with coaction
$X\xra{\rho_X}C\otimes X$ (resp. $X\xra{\rho_X}X\otimes C$), we write
$\rho_X(x)=x_{(-1)}\otimes x_{(0)}$ (resp.  $\rho_X(x)=x_{(0)}\otimes
x_{(1)}$) to denote its structure map where a summation is understood.

\section{Module coalgebras and equivariant comodules}

Let $B$ be a bialgebra. Recall that a coalgebra $C$ is called a left
$B$--module coalgebra if $C$ is a left $B$-module and its
comultiplication $C\xra{\Delta}C\otimes C$ is a $B$--module morphism.
Explicitly
\begin{align*}
  (b(c))_{(1)}\otimes (b(c))_{(2)}
  = b_{(1)}(c_{(1)})\otimes b_{(2)}(c_{(2)})
\end{align*}
for any $c\in C$ and $b\in B$.

  Let $C$ be a coalgebra, $X\xra{\rho_X}X\otimes C$ be a right
  $C$--comodule and $Y\xra{\rho_Y}C\otimes Y$ be a left $C$--comodule.
  Recall that the cotensor product $X\fotimes{C}Y$ is defined as the
  kernel
  \begin{align*}
    X\fotimes{C}Y := ker\left(X\otimes Y
      \xra{(\rho_X\otimes id_Y)-(id_X\otimes\rho_Y)}X\otimes C\otimes Y\right)
  \end{align*}
  Explicitly, $X\fotimes{C}Y$ is the $k$--submodule of $X\otimes Y$
  generated over $k$ by terms of the form $\sum_i (x^i\otimes y^i)$
  which satisfy
  \begin{align*}
    \sum_i x^i_{(0)}\otimes x^i_{(1)}\otimes y
      = \sum_i x^i\otimes y^i_{(-1)}\otimes y^i_{(0)}
  \end{align*}
  The derived functors of the cotensor product are denoted by ${\rm
  Cotor}^C_*(\cdot,\cdot)$.

\begin{defn}
  Let $C$ be a left $B$--module coalgebra.  A $k$--module $M$ is
  called a $B$--equivariant left $C$--comodule if: (i) $M$ is a left
  $B$--module, (ii) $M$ is a left $C$--comodule, and (iii) the
  $C$--coaction $M\xra{\rho_M}C\otimes M$ is a $B$--module morphism,
  in other words
  \begin{align*}
    (b(m))_{(-1)}\otimes (b(m))_{(0)} = b_{(1)}(m_{(-1)})\otimes
     b_{(2)}(m_{(0)})
  \end{align*}
  for all $b\in B$ and $m\in M$.
\end{defn}

\begin{defn}
  Let $C$ be a coalgebra.  We define a differential graded $k$--module
  ${\rm CB}^{\rm bar}_*(C)$ by letting ${\rm CB}^{\rm
  bar}_n(C)=C^{\otimes n+2}$ and defining a pre-simplicial structure
  \begin{align*}
    \partial_j(c^0\otimes\cdots\otimes c^{n+1})
    = c^0\otimes \cdots\otimes c^j_{(1)}\otimes c^j_{(2)}\otimes\cdots \otimes c^{n+1}
  \end{align*}
  for any $n\geq 0$, $0\leq j\leq n+1$ and $(c^0\otimes\cdots\otimes
  c^{n+1})$ from ${\rm CB}^{\rm bar}_n(C)$.  The corresponding
  differential is defined as $d^{\rm CB}_n
  =\sum_{j=0}^{n+1}(-1)^j\partial_j$. The differential graded
  $k$--module ${\rm CB}^{\rm bar}_*(C)$ is called the bar resolution
  of $C$.  Note that ${\rm CB}^{\rm bar}_*(C)$ need not be a
  resolution for a non-counital coalgebra $C$.
\end{defn}

\begin{prop}\label{EquivariantBarComplex}
  Let $C$ be a left $B$--module coalgebra.  Then the bar resolution
  ${\rm CB}^{\rm bar}_*(C)$ is a differential graded $B$--equivariant
  $C$--bicomodule.
\end{prop}

\begin{proof}
  The left and right $C$--comodule structures are defined as
  \begin{align*}
    \rho_L(c^0\otimes\cdots\otimes c^{n+1})
    = & c^0_{(1)}\otimes (c^0_{(2)}\otimes\cdots\otimes c^{n+1})\\
    \rho_R(c^0\otimes\cdots\otimes c^{n+1})
    = & (c^0\otimes\cdots\otimes c^{n+1}_{(1)})\otimes c^{n+1}_{(2)}
  \end{align*}
  and the $B$--module structure is diagonal
  \begin{align*}
    L_b(c^0\otimes\cdots\otimes c^{n+1})
    = b_{(1)}(c^0)\otimes\cdots\otimes b_{(n+1)}(c^{n+1})
  \end{align*}
  defined for any $b\in B$, $n\geq 0$ and $c^0\otimes\cdots\otimes
  c^{n+1}$ from ${\rm CB}^{\rm bar}_n(C)$.  As for the $C$--coaction
  being $B$--equivariant we check
  \begin{align*}
    L_b\rho_L(c^0\otimes\cdots\otimes c^{n+1})
    = & b_{(1)}(c^0_{(1)})\otimes b_{(2)}(c^0_{(2)})\otimes
        b_{(3)}(c^1)\otimes\cdots\otimes b_{(n+3)}(c^{n+1})\\
    = & b_{(1)}(c^0_{(1)})\otimes L_{b_{(2)}}(c^0_{(2)}\otimes\cdots\otimes c^{n+1})\\
    = & \rho_L L_b(c^0\otimes\cdots\otimes c^{n+1})
  \end{align*}
  for any $b\in B$, $n\geq 0$ and $(c^0\otimes\cdots\otimes c^{n+1})$
  from ${\rm CB}^{\rm bar}_n(C)$.  The proof for the right coaction is
  similar.  So far we only showed that ${\rm CB}^{\rm bar}_*(C)$ is a
  graded $B$--equivariant $C$--bicomodule.  We must show that the
  pre-simplicial structure is compatible with the $C$--bicomodule
  structure and the $B$--module structures.  We start with the
  $B$--module structure
  \begin{align*}
    L_b\partial_j & (c^0\otimes\cdots\otimes c^{n+1})\\
    = & \begin{cases}
        b_{(1)}(c^0_{(1)})\otimes b_{(2)}(c^0_{(2)})\otimes
        b_{(3)}(c^1)\otimes\cdots\otimes b_{(n+3)}(c^{n+1})
      & \text{ if } j=0\\
        b_{(1)}(c^0)\otimes\cdots\otimes
         b_{(j+1)}(c^j_{(1)})\otimes b_{(j+2)}(c^j_{(2)})\otimes
         \cdots\otimes b_{(n+3)}(c^{n+1})
      & \text{ if } 1\leq j\leq n\\
        b_{(1)}(c^0)\otimes\cdots\otimes b_{(n+1)}(c^n)\otimes
        b_{(n+2)}(c^{n+1}_{(1)})\otimes b_{(n+3)}(c^{n+1}_{(2)})
      & \text{ if } j=n+1
    \end{cases}\\
    = & \partial_jL_b(c^0\otimes\cdots\otimes c^{n+1})
  \end{align*}
  As for the $C$--comodule structure
  \begin{align*}
    (id_C\otimes\partial_j)\rho_L(c^0\otimes\cdots\otimes c^{n+1})
    = & \begin{cases}
          c^0_{(1)}\otimes c^0_{(2)}\otimes c^0_{(3)}\otimes
            c^1\otimes\cdots\otimes c^{n+1}
            & \text{ if } j=0\\
          c^0_{(1)}\otimes c^0_{(2)}\otimes\cdots\otimes c^j_{(1)}\otimes
            c^j_{(2)}\otimes\cdots
            & \text{ if } 1\leq j\leq n+1
    \end{cases}\\
    = & \rho_L\partial_j(c^0\otimes\cdots\otimes c^{n+1})
  \end{align*}
  The proof for the right coaction is similar.
\end{proof}

\begin{lem}\label{CotorProduct}
  Let $C$ be a $B$--module coalgebra, $X$ be a $B$--equivariant right
  $C$--comodule and $Y$ be a $B$--equivariant left $C$--comodule.
  Then $X\fotimes{C}Y$ is a $B$--module.
\end{lem}

\begin{proof}
  The left $B$ action on $X\otimes Y$ is diagonal.  We check that
  $X\fotimes{C}Y$ is a $B$--submodule.  Take $\Psi$ from
  $X\fotimes{C}Y$ and $b\in B$.  Then
  \begin{align*}
    (\rho_X\otimes id_Y)L_b(\Psi)
    = & L_b(\rho_X\otimes id_Y)(\Psi)
    = L_b(id_X\otimes\rho_Y)(\Psi)
    = (id_X\otimes\rho_Y)L_b(\Psi)
  \end{align*}
  proving what we wanted.
\end{proof}

\begin{defn}
  An object $X$ in the category of right (resp. left) $C$--comodules
  is called co-projective if the functor
  $\text{Hom}_{\text{Comod--}C}(X,\cdot)$
  (resp. $\text{Hom}_{C\text{--Comod}}(X,\cdot)$) is an exact functor
  from the category of right (resp. left) $C$--comodules into the
  category of $k$--modules.
\end{defn}

\begin{prop}
  Let $C$ be a counital $B$--module coalgebra and $X$ and $Y$ be as
  before.  Then the cotorsion groups ${\rm Cotor}^C_*(X,Y)$ are
  canonically $B$--modules.
\end{prop}

\begin{proof}
  The existence of right derived functors of the functors
  $X\fotimes{C}\cdot$ and $\cdot\fotimes{C}Y$ hinge on the fact that
  we have enough co-projectives, i.e. each object in the category of
  left (resp. right) $C$--comodules admits a co-projective resolution.
  Since $C$ is counital one has $X\fotimes{C}C\cong X$ and
  $C\fotimes{C}Y\cong Y$ as $B$--modules meaning $C$ is co-projective
  in the category of $B$--equivariant $C$--comodules both left and
  right.  Since $k$ is a field, ${\rm CB}^{\rm bar}_*(C)\fotimes{C}Y$
  and $X\fotimes{C}{\rm CB}^{\rm bar}_*(C)$ will provide these
  resolutions for $Y$ and $X$ respectively.  Then the derived functors
  of the cotorsion product exist and one can immediately see that
  ${\rm Cotor}^C_*(X,Y)$ are naturally $B$--modules after
  Lemma~\ref{CotorProduct} and
  Proposition~\ref{EquivariantBarComplex}.
\end{proof}

\begin{defn}
  The bar complex ${\rm CB}_*(C)$ associated with a coalgebra $C$ is
  defined as ${\rm CB}_n(C)={\rm CB}^{\rm bar}_{n-1}(C)$ for any
  $n\geq 1$ and ${\rm CB}_0(C)=C$.  The differentials too are the same
  for $n\geq 1$ and we let $d^{\rm CB}_0=\Delta_C$.
\end{defn}

\begin{defn}
  A coalgebra $C$ is called H--counital if ${\rm CB}_*(C)$ is acyclic.
  This is equivalent to saying the bar resolution ${\rm CB}^{\rm
  bar}_*(C)$ is quasi-isomorphic to the differential graded
  $k$--module $C$ concentrated at degree $0$.
\end{defn}

\section{Twisted Cartier--Hochschild homology}

  Let $C$ be a coalgebra.  Define a coalgebra $C^e=C\otimes C^{op}$
  with the diagonal comultiplication
  \begin{align*}
    \Delta(c\otimes c')
      = (c_{(1)}\otimes c'_{(2)})\otimes (c_{(2)}\otimes c'_{(1)})
  \end{align*}
  for any $c\otimes c'$ from $C^e$.

  It is easy to see that the categories of $C$--bicomodules and right
  $C^e$--comodules are isomorphic. Given any $C$--bicomodule $M$ the
  corresponding right $C^e$--comodule is again $M$ with comodule
  structure defined as
  \begin{align*}
    \rho^e_M(m) = m_{(0)}\otimes (m_{(1)}\otimes m_{(-1)})
  \end{align*}
  for any $m\in M$.  

  Assume $C$ is a $B$--module coalgebra and $X$ and $Y$ are left
  $B$--equivariant $C$--bicomodules.  Note that unless $B$ is
  cocommutative, the cotensor product $X\fotimes{C^e}Y$ is not a
  $B$--module under the diagonal $B$ action. 

\begin{defn}
  Let $C$ be a coalgebra and $M$ be a $C$--bicomodule.  Define a
  differential graded $k$--module ${\rm CH}_*(C,M)$ by letting ${\rm
  CH}_n(C,M)=M\otimes C^{\otimes n}$.  We define a pre-simplicial
  structure
  \begin{align*}
    \partial_j(m\otimes c^1\otimes\cdots\otimes c^n)
    = \begin{cases}
      m_{(0)}\otimes m_{(1)}\otimes c^1\otimes\cdots\otimes c^n
         & \text{ if } j=0\\
      m\otimes c^1\otimes\cdots\otimes c^j_{(1)}\otimes c^j_{(2)}
         \otimes\cdots\otimes c^n
         & \text{ if } 1\leq j\leq n\\
      m_{(0)}\otimes c^1\otimes\cdots\otimes c^n\otimes m_{(-1)}
         & \text{ if } j=n+1
      \end{cases}
  \end{align*}
  and a differential $d^{\rm CH}_n=\sum_{j=0}^{n+1}(-1)^j\partial_j$.
\end{defn}

  Note that ${\rm CH}_*(C,M)$ is a
  differential graded $k$--module and a graded $B$--module but not a
  differential graded $B$--module unless $B$ is cocommutative.

\begin{thm}[\cite{Doi:HomologicalCoalgebra}]\label{Doi}
  Let $C$ be a counital coalgebra and $M$ be a $C$--bicomodule.  Then
  there is an isomorphism of differential graded modules
  ${\rm CH}_*(C,M)\cong M\fotimes{C^e}{\rm CB}^{\rm bar}_*(C)$
\end{thm}

\begin{cor}\label{Cofibration}
  Let $C$ be a counital coalgebra and
  $0\xra{}N\xra{}M\xra{}M/N\xra{}0$ be a short exact sequence of
  $C$--bicomodules.  Then we have a short exact sequence of
  differential graded $k$--modules of the form
  \begin{align*}
    0\xra{}{\rm CH}_*(C,N)\xra{}{\rm CH}_*(C,M)\xra{}{\rm CH}_*(C,M/N)\xra{}0
  \end{align*}
\end{cor}

\begin{proof}
  Since ${\rm CB}^{\rm bar}_*(C)$ consists of co-projective
  $C^e$--comodules we get a short exact sequence of the form
  \begin{align*}
    0\xra{}N\fotimes{C^e}{\rm CB}^{\rm bar}_*(C)
     \xra{}M\fotimes{C^e}{\rm CB}^{\rm bar}_*(C)
     \xra{}M/N\fotimes{C^e}{\rm CB}^{\rm bar}_*(C)\xra{}0
  \end{align*}
  The result follows after Theorem~\ref{Doi}.
\end{proof}

\begin{rem}
  One can give a direct proof of Corollary~\ref{Cofibration} without
  appealing to $\cdot\fotimes{C^e}{\rm CB}^{\rm bar}_*(C)$.  Indeed,
  for a direct proof one need not $C$ to be counital.
\end{rem}

\begin{defn}
  Let $X$ be a left $B$--comodule and $M$ be a $C$--bicomodule.
  Define $X\rtimes M$ as the $C^e$--comodule with $X\rtimes M =
  X\otimes M$ on the $k$--module level and we let
  \begin{align*}
    \rho^e(x\otimes m)
      & = (x_{(0)}\otimes m_{(0)})\otimes (m_{(1)}\otimes x_{(-1)} m_{(-1)})
  \end{align*}
  Note that
  \begin{align*}
    (\rho^e\otimes id_{C^e})\rho^e(x\otimes M)
    & = (x_{(0)(0)}\otimes m_{(0)(0)})\otimes (m_{(0)(1)}\otimes x_{(0)(-1)} m_{(0)(-1)})
         \otimes (m_{(1)}\otimes x_{(-1)} m_{(-11)}) \\
    & = (x_{(0)}\otimes m_{(0)})\otimes (m_{(1)}\otimes x_{(-1)}m_{(-1)})
         \otimes (m_{(2)}\otimes x_{(-2)} m_{(-2)})\\
    & = (x_{(0)}\otimes m_{(0)})\otimes (m_{(1)(1)}\otimes x_{(-1)(2)}m_{(-1)(2)})
         \otimes (m_{(1)(2)}\otimes x_{(-1)(1)} m_{(-1)(1)})\\
    & = (id_{X\rtimes C}\otimes \Delta)\rho^e(x\otimes m)
  \end{align*}
  for any $(x\otimes m)$ from $X\rtimes M$.
\end{defn}

\begin{defn}\label{CartierHochschild}
  Assume $X$ is a $B$--module/comodule and $M$ is a $C$--bicomodule.
  We define the Cartier-Hochschild complex of  $C$ with
  coefficients in $M$ twisted by $X$ as
  \begin{align*}
    {\rm CH}^{\rm tw}_*(C,M;X)
    := {\rm CH}_*(C,X\rtimes M)
    := (X\rtimes M)\fotimes{C^e}{\rm CB}^{\rm bar}_*(C)
  \end{align*}
  In the special case where $M=C$ we will write
  ${\rm CH}^{\rm tw}_*(C;X):={\rm CH}^{\rm tw}_*(C,C;X)$.
\end{defn}

\begin{prop}\label{GradedModule}
  Let $C$ be a counital $B$--module coalgebra and $M$ be a
  $B$--equivariant $C$--bicomodule.  There is a graded $B$--module
  structure on ${\rm CH}^{\rm tw}_*(C,M;X)$.
\end{prop}

\begin{proof}
  We define a $B$--module structure on ${\rm CH}^{\rm tw}_*(C,M;X)$ as
  follows: we let
  \begin{align*}
    L_b(x\otimes m\otimes c^1\otimes\cdots\otimes c^n)
    = b_{(n+2)}x\otimes b_{(1)}m\otimes b_{(2)}c^1\otimes\cdots\otimes b_{(n+1)}c^n
  \end{align*}
  for any $n\geq 0$, $x\otimes m\otimes c^1\otimes\cdots\otimes c^n$
  from ${\rm CH}^{\rm tw}_n(C,M;X)$ and $b\in B$.  One can easily
  check that the action is coassociative and $[L_b,\partial_j]=0$ for
  any $b\in B$, $0\leq n$ and $0\leq j\leq n$.  Note that, unless we
  have a prescribed interaction between the $B$--action and
  $B$--coaction on $X$, we can not infer anything about
  $[L_b,\partial_{n+1}]$.
\end{proof}

Let us recall the following definition from
\cite{Khalkhali:SaYDModules}.
\begin{defn}
  Let $B$ be a Hopf algebra with an invertible antipode.  Assume $X$
  is a left $B$--module  and a left
  $B$--comodule.  $X$ is called  an
  anti-Yetter-Drinfeld (aYD) module if one has
  \begin{align*}
    (bx)_{(-1)}\otimes (bx)_{(0)}
     = b_{(1)}x_{(-1)}S^{-1}(b_{(3)})\otimes b_{(2)}x_{(0)}
  \end{align*}
  for any $x\in X$ and $b\in B$.
\end{defn}

\begin{thm}\label{aYD}
  Let $B$ be a Hopf algebra with an invertible antipode and $X$ be an
  aYD module.  Assume $C$ is a $B$--module coalgebra and $M$ is a
  $B$--equivariant $C$--bicomodule.  Then $k\eotimes{B}{\rm CH}^{\rm
  tw}_*(C,M;X)$ is a differential graded $k$--module with the induced
  differentials $id\eotimes{B} d^{\rm CH}_*$.
\end{thm}

\begin{proof}
  Consider the graded $k$--submodule $I_*$ of ${\rm CH}^{\rm
  tw}_*(C,M;X)$ generated by the commutators of the form $[L_b,d^{\rm
  CH}_n](\Psi)$ where $n\geq 0$, $\Psi\in{\rm CH}^{\rm tw}_*(C,M;X)$
  and $b\in B$.  Let us show that $I_*$ is a differential graded
  $k$--submodule and also a graded $B$--submodule of ${\rm CH}^{\rm
  tw}_*(C,M;X)$.  Start with
  \begin{align*}
    d^{\rm CH}_{n+1}[L_b,d^{\rm CH}_n]
    = - [L_b,d^{\rm CH}_{n+1}]d^{\rm CH}_n + [L_b,d^{\rm CH}_{n+1}d^{\rm CH}_n]
    = - [L_b,d^{\rm CH}_{n+1}]d^{\rm CH}_n
  \end{align*}
  which proves $I_*$ is a differential graded $k$--submodule.  Then
  \begin{align*}
    L_y[L_b,d^{\rm CH}_n]
    = -[L_y,d^{\rm CH}_n]L_b + [L_{by},d^{\rm CH}_n]
  \end{align*}
  which proves $I_*$ is a graded $B$--submodule.  By design ${\rm
  CH}^{\rm tw}_*(C,M;X)/I_*$ is a differential graded $B$--module.
  Therefore $k\eotimes{B}\left({\rm CH}^{\rm tw}_*(C,M;X)/I_*\right)$
  is a differential graded $k$--module.  The rest of the proof
  concentrates on proving $k\eotimes{B}I_*=0$ which would imply that
  the graded $k$--modules $k\eotimes{B}\left({\rm CH}^{\rm
  tw}_*(C,M;X)/I_*\right)$ and $k\eotimes{B}{\rm CH}^{\rm
  tw}_*(C,M,X)$ are isomorphic.  Note the remark we made in the proof
  of Proposition~\ref{GradedModule} and observe that $[L_b,d^{\rm
  CH}_n]=(-1)^{n+1}[L_b,\partial_{n+1}]$ for any $n\geq 0$.  Take an
  arbitrary element $\Phi$ from $k\eotimes{B}I_*$ of the form $\sum_i
  1\eotimes{B}[L_{b^i},\partial_{n+1}](x^i\otimes m^i\otimes
  c^{1,i}\otimes\cdots\otimes c^{n,i})$ and consider the equalities
  \begin{align*}
    1\eotimes{B} & \partial_{n+1}L_{b^i}(x^i\otimes m^i\otimes
        c^{1,i}\otimes\cdots\otimes c^{n,i})
    =   1\eotimes{B}\partial_{n+1}(b^i_{(n+2)}x^i\otimes b^i_{(1)}m^i
         \otimes b^i_{(2)}c^{1,i}\otimes\cdots\otimes b^i_{(n+1)}c^{n,i})\\
    = & 1\eotimes{B}b^i_{(n+4)}x^i_{(0)}\otimes b^i_{(2)}m^i_{(0)}\otimes
         \otimes b^i_{(3)}c^{1,i}\otimes\cdots\otimes b^i_{(n+2)}c^{n,i}
         \otimes b^i_{(n+3)}x_{(-1)}S^{-1}(b^i_{(n+5)})b^i_{(1)}m_{(-1)}\\
    = & 1\eotimes{B}L_{b^i_{(2)}}(x^i_{(0)}\otimes m^i_{(0)}\otimes c^{1,i}
         \otimes\cdots\otimes c^{n,i}\otimes x_{(-1)}S^{-1}(b^i_{(3)})b^i_{(1)}m_{(-1)})\\
    = & \varepsilon(b^i)\eotimes{B} (x^i_{(0)}\otimes m^i_{(0)}\otimes c^{1,i}
         \otimes\cdots\otimes c^{n,i}\otimes x_{(-1)}m_{(-1)})\\
    = & \varepsilon(b^i)\eotimes{B}\partial_{n+1}(x^i\otimes m^i\otimes c^{1,i}
         \otimes\cdots\otimes c^{n,i})
  \end{align*}
  which means $\Phi=0$ as we wanted to show.
\end{proof}

\begin{thm}\label{Hochschild0}
  Assume $B$ is a Hopf algebra with an invertible antipode and $X$ is
  an aYD module.  Let $C$ be a counital $B$--module coalgebra and
  $0\xra{}N\xra{}M\xra{}M/N\xra{}0$ be a short exact sequence of
  $B$--equivariant $C$--bicomodules.  Assume $N\xra{}M$ is a split
  injective morphism of $B$--modules.  Then there is a short exact
  sequence of differential graded $k$--modules of the form
  \begin{align}\label{HochschildCofibration}
    0\xra{}k\eotimes{B}{\rm CH}^{\rm tw}_*(C,N;X)
     \xra{}k\eotimes{B}{\rm CH}^{\rm tw}_*(C,M;X)
     \xra{}k\eotimes{B}{\rm CH}^{\rm tw}_*(C,M/N;X)\xra{}0
  \end{align}
\end{thm}

\begin{proof}
  Corollary~\ref{Cofibration} implies we have a short exact sequence
  of differential graded $k$--modules and graded $B$--modules
  \begin{align*}
    0\xra{}{\rm CH}^{\rm tw}_*(C,N;X)
     \xra{}{\rm CH}^{\rm tw}_*(C,M;X)
     \xra{}{\rm CH}^{\rm tw}_*(C,M/N;X)\xra{}0
  \end{align*}
  Upon tensoring with $k$ over $B$ we get the short exact sequence in
  Equation~\ref{HochschildCofibration} as a short exact sequence of
  graded $k$--modules.  Theorem~\ref{aYD} tells us that each of the
  terms in the short exact sequence is a differential graded
  $k$--module.  Since the differentials are induced the short exact
  sequence in Equation~\ref{HochschildCofibration} is a short exact
  sequence of differential graded $k$--modules.
\end{proof}

\section{Excision in Hopf cyclic homology}

  In this section we assume $B$ is a Hopf algebra with an invertible
  antipode and $X$ is an aYD module unless otherwise stated.

\begin{lem}\label{MacLane}
  Let $P_1,\ldots,P_n$ be a finite collection of projective
  $B$--modules.  Then $P_1\otimes\cdots\otimes P_n$ is projective with
  diagonal $B$--action.
\end{lem}

\begin{proof}
  Since every projective module is a direct summand of a free module,
  and direct sum of projective modules is projective, the statement
  reduces to proving $B^{\otimes n}$ is projective.  Consider the
  endomorphism of $k$--modules $B^{\otimes n}\xra{\Gamma}B^{\otimes
  n}$ defined by
  \begin{align*}
    \Gamma_n(b^1\otimes\cdots\otimes b^n)
      & = b^1_{(1)}\otimes b^1_{(2)}b^2_{(1)}\otimes\cdots\otimes
          b^1_{(n)}b^2_{(n-1)}\cdots b^{n-1}_{(2)}b^n
  \end{align*}
  Note that $\Gamma_n(bb^1\otimes\cdots\otimes b^n)
  =L_b\Gamma_n(b^1\otimes\cdots\otimes b^n)$ for any $b\in B$ and
  $(b^1\otimes\cdots\otimes b^n)$ from $B^{\otimes n}$.  $\Gamma$ has
  an inverse given by
  \begin{align*}
    \Gamma_n^{-1}(b^1\otimes\cdots\otimes b^n)
     = b^1_{(1)}\otimes S(b^1_{(2)})b^2_{(1)}\otimes\cdots\otimes
         S(b^{n-1}_{(2)})b^n
  \end{align*}
  which means $B^{\otimes n}$ with diagonal $B$--action is isomorphic
  to the $B$--module $B^{\otimes n}$ with $B$ acting on the first
  tensor component.  Since $k$ is a field, the second $B$--module is
  $B$--free.  The result follows.
\end{proof}

\begin{lem}\label{Trivialization}
  Given a left $B$--module $U$ define right $B$--module $U^{op}$ by
  letting $u\cdot b := S^{-1}(b) u$ for any $u\in U$ and $b\in B$.
  Then given any other left $B$--module $V$ there is an isomorphism of
  $k$--modules
  \begin{align*}
    k\eotimes{B}(U\otimes V)\cong U^{op}\eotimes{B} V
  \end{align*}
\end{lem}

\begin{proof}
  On one direction we define
  $U^{op}\eotimes{B}V\xra{\phi}k\eotimes{B}(U\otimes V)$ by sending
  each $u^{op}\eotimes{B} v$ in $U^{op}\eotimes{B}V$ to
  $1\eotimes{B}(u\otimes v)$.  We must show that $\phi$ is
  well-defined.  For that, for any $(u^{op}b\eotimes{B} v)$ from
  $U^{op}\eotimes{B}V$ we consider
  \begin{align*}
    \phi(u^{op}b\eotimes{B} v)
     = 1\eotimes{B}(S^{-1}(b)u\otimes v)
     = 1\eotimes{B}S^{-1}(b_{(2)})(u\otimes b_{(1)}v)
     = 1\eotimes{B}(u\otimes bv)
     = \phi(u^{op}\eotimes{B}bv)
  \end{align*}
  which proves $\phi$ is well-defined.  Conversely take
  $1\eotimes{B}(u\otimes v)$ from $k\eotimes{B}(U\otimes V)$ and let
  $\psi(1\eotimes{B}(u\otimes v))$ be $(u^{op}\eotimes{B}v)$.  In
  order $\psi$ be well-defined we must have
  \begin{align*}
    \psi(1\eotimes{B}x(u\otimes v))
     = (u^{op}S(x_{(1)})\eotimes{B}x_{(2)}v)
     = \varepsilon(x)(u^{op}\eotimes{B}v)
     = \psi(\varepsilon(x)\eotimes{B}(u\otimes v))
  \end{align*}
  we we wanted to show.  One can easily see that $\psi$ and $\phi$ are
  mutual inverses.
\end{proof}

\begin{lem}\label{FiniteDimension}
  Let $\C{A}_*$ be an acyclic differential graded $B$--module where
  each $\C{A}_n$ is $B$--projective.  Let $U$ be a $B$--module with
  finite projective dimension and let $V$ be a projective $B$--module.
  Then $k\eotimes{B}(U\otimes\C{A}_*\otimes V)$ is acyclic.
\end{lem}

\begin{proof}
  We will show that $U^{op}\eotimes{B}(\C{A}_*\otimes V)$ is acyclic
  which after using Lemma~\ref{Trivialization} would imply
  $k\eotimes{B}(U\otimes\C{A}_*\otimes V)$ is acyclic.  Note that
  $\C{A}_*\otimes V$ is an acyclic differential graded $B$--module
  which consists of projective $B$--modules following
  Lemma~\ref{MacLane}.  Since each $\C{A}_n\otimes V$ is projective,
  the differential graded $B$--module $\C{A}_n\otimes V$ is its own
  Cartan--Eilenberg resolution.  Then the hyper Tor products are given
  by
  \begin{align*}
    {\rm\bf Tor}^B_*(U^{op},\C{A}_*\otimes V)
     = H_*(U^{op}\eotimes{B}(\C{A}_*\otimes V))
  \end{align*}
  Since we actually compute the hyper Tor functors, one can compute
  the homology of $U^{op}\eotimes{B}(\C{A}_*\otimes V)$ by taking a
  $B$--projective resolution $(\C{U}_*,d^\C{U}_*)$ of $U^{op}$ and
  then compute the homology of the double complex
  $\C{U}_*\eotimes{B}(\C{A}_*\otimes V)$.  Since we did not assume
  $\C{A}_*$ is bounded, there is a problem about the convergence of
  the canonical spectral sequence
  \begin{align*}
    E^2_{pq} = H_p(\C{U}_*\eotimes{B}(H_q(\C{A}_*,d^\C{A}_*)\otimes V), d^\C{U}_*)
  \end{align*}
  coming from this double complex.  However, since we assumed $U$ is
  of finite projective dimension, the double complex
  $\C{U}_*\eotimes{B}(\C{A}_*\otimes V)$ is concentrated on an
  infinite {\rm horizontal} strip which forces the spectral sequence
  to converge.  Since $\C{A}_*$ is acyclic, the spectral sequence is
  all $0$.  Thus $U^{op}\eotimes{B}(\C{A}_*\otimes V)$ is acyclic.
\end{proof}

\begin{thm}\label{Hochschild1}
  Let $0\xra{}K\xra{}C\xra{}C/K\xra{}0$ be a short exact sequence of
  $B$--module coalgebras.  Assume that both $C$ and $C/K$ are
  $B$--projective and $X$ has finite projective dimension. Assume also
  that $C/K$ is H--counital.  Then the canonical epimorphism
  $C\xra{}C/K$ of $B$--module coalgebras induces a weak equivalence of
  differential graded $k$--modules of the form $k\eotimes{B}{\rm
  CH}^{\rm tw}_*(C,C/K;X)\xra{} k\eotimes{B}{\rm CH}^{\rm
  tw}_*(C/K;X)$.
\end{thm}

\begin{proof}
  First note that since $C/K$ is projective the short exact sequence
  of $B$--module coalgebras $0\xra{}K\xra{}C\xra{}C/K\xra{}0$ is a
  split short exact sequence of $B$--modules.  Define a sequence
  $\{G^p_*\}_{p\geq 0}$ of differential graded $k$--modules by letting
  \begin{align*}
    G^p_* = \begin{cases}
            X\otimes (C/K)\otimes C^{\otimes n} & \text{ if } n\leq p\\
            X\otimes (C/K)^{\otimes n-p+1}\otimes C^{\otimes p} & \text{ if } n>p
        \end{cases}
  \end{align*}
  The pre-simplicial structure for $n< p$ comes from ${\rm CH}^{\rm
  tw}_*(C,C/K;X)$.  For $n\geq p$
  \begin{align*}
    \partial_j
     = \begin{cases}
       (id_X\otimes id_j\otimes\Delta_{C/K}\otimes id_{n-j})
             & \text{ if } 0\leq j\leq n-p\\
       (id_X\otimes id_{n-p+1}\otimes\pi_{C/K}\otimes id_p)
       (id_X\otimes id_j\otimes\Delta_C\otimes id_{n-j})
             & \text{ if } n-p+1\leq j\leq n\\
       (id_X\otimes id_{n-p+1}\otimes\pi_{C/K}\otimes id_p)
       (id_X\otimes\tau_{n+1})(\rho_{C/K}\otimes id_n)
             & \text{ if } j=n+1
       \end{cases}
  \end{align*}
  where $C\xra{\pi_{C/K}}C/K$ is the canonical epimorphism and
  $\tau_{n+1}$ is the twisted cyclic permutation of length $n+2$ and
  $C/K\xra{\rho_{C/K}}C\otimes C/K$ is the left $C$--comodule
  structure on $C/K$.  There are epimorphisms of differential graded
  $k$--modules $G^{p+1}_*\xra{\pi^p_*}G^p_*$ defined for any $p\geq 0$
  where
  \begin{align*}
    \pi^p_n = \begin{cases}
              (id_X\otimes id_{n+1})  & \text{ if } n\leq p\\
              (id_X\otimes id_{n-p}\otimes\pi_{C/K}\otimes id_p) & \text{ if } n>p
          \end{cases}
  \end{align*}
  with kernel
  \begin{align*}
    ker(\pi^p_*) =
       \begin{cases}
     0 & \text{ if } n\leq p\\
         X\otimes (C/K)^{\otimes n-p}\otimes K\otimes C^{\otimes p}
           & \text{ if } n>p
       \end{cases}
  \end{align*}
  which is isomorphic as graded $B$--modules to $X\otimes {\rm
  CB}_*(C/K)[-p]\otimes K\otimes C^{\otimes p}$.  Now, applying the
  functor $k\eotimes{B}\cdot$ we get a short exact sequence
  \begin{align*}
    0\xra{}k\eotimes{B}ker(\pi^p_*)\xra{}k\eotimes{B}G^{p+1}_*
     \xra{id\eotimes{B}\pi^p_*}k\eotimes{B}G^p_*\xra{}0
  \end{align*}
  of graded $k$--modules by using the fact that the inclusion
  $K\xra{}C$ is $B$--split.  Theorem~\ref{aYD} implies
  $k\eotimes{B}G^p_*$ is a differential graded $k$--module with the
  induced differentials.  Since $C/K$ is H--counital
  Lemma~\ref{MacLane} and Lemma~\ref{FiniteDimension} implies
  $k\eotimes{B}ker(\pi^p_*)$ is acyclic. Now observe that
  $k\eotimes{B}G^0_*=k\eotimes{B}{\rm CH}^{\rm tw}(C/K;X)$ and
  $k\eotimes{B}G^p_n = k\eotimes{B}{\rm CH}^{\rm tw}_n(C,C/K;X)$ for
  $n>p$.  The result follows.
\end{proof}

\begin{prop}\label{ContractibleCotor}
  Let $0\xra{}K\xra{}C\xra{}C/K\xra{}0$ be a short exact sequence of
  $B$--module coalgebras.  Assume that $C/K$ is H--counital.  We
  define a differential graded $B$--module ${\rm C}_*(K,C/K)$ by
  letting ${\rm C}_n(K,C/K)=K\otimes C^{\otimes n}\otimes C/K$ then
  defining the pre-simplicial structure as
  \begin{align*}
    \partial_j = \begin{cases}
                 \Delta_K\otimes id_{n+1} & \text{ if } j=0\\
                 id_{j+1}\otimes\Delta_C\otimes id_{n-j+1} & \text{ if } 0<j<n+1\\
                 id_{n+1}\otimes\rho_{C/K} & \text{ if } j=n+1
             \end{cases}
  \end{align*}
  where $\Delta_C$ and $\Delta_K$ are the comultiplication on $C$ and
  $K$ respectively, and $C/K\xra{\rho_{C/K}}C\otimes C/K$ is the
  $C$--comodule structure on $C/K$.  Then ${\rm C}_*(K,C/K)$ is an
  acyclic differential graded $B$--module.
\end{prop}

\begin{proof}
  The fact that ${\rm C}_*(K,C/K)$ is a differential graded
  $B$--module is easy to see.  Consider the sequence of differential
  graded $B$--modules
  \begin{align*}
    J^p_* = \begin{cases}
      K\otimes C^{\otimes n}\otimes C/K & \text{ if } n\leq p\\
      K\otimes C^{\otimes p}\otimes(C/K)^{\otimes n-p+1}
                                        & \text{ if } n>p
            \end{cases}
  \end{align*}
  The differentials on $J^p_*$ descend from ${\rm C}_*(K,C/K)$ for
  $n<p$.  For $n\geq p$ we define a pre-simplicial structure
  \begin{align*}
    \partial_j =
      \begin{cases}
        (id_{p+1}\otimes\pi_{C/K}\otimes id_{n-p})(\Delta_K\otimes id_{n+1})
                               & \text{ if } j=0\\
        (id_{p+1}\otimes\pi_{C/K}\otimes id_{n-p})(id_j\otimes\Delta_C\otimes id_{n-j+1})
                                   & \text{ if } 1\leq j\leq p\\
        (id_j\otimes\Delta_{C/K}\otimes id_{n-j+1})
                                   & \text{ if } p<j\leq n+1
      \end{cases}
  \end{align*}
  There are epimorphisms of differential graded $B$--modules
  $J^{p+1}_*\xra{\pi^p_*}J^p_*$ defined by
  \begin{align*}
    \pi^p_n =  \begin{cases}
               id_{n+2}  & \text{ if } n\leq p\\
           id_{p+1}\otimes\pi_{C/K}\otimes id_{n-p} & \text{ if } n>p
               \end{cases}
  \end{align*}
  with kernels
  \begin{align*}
    ker(\pi^p_*)
    = \begin{cases}
      0 & \text{ if } n\leq p\\
      K\otimes C^{\otimes p}\otimes K\otimes (C/K)^{\otimes n-p} & \text{ if } n>p
      \end{cases}
  \end{align*}
  which is isomorphic to the differential graded $B$--module $K\otimes
  C^{\otimes p}\otimes K\otimes {\rm CB}_*(C/K)[-p]$.  Since $C/K$ is
  H--counital, we see that each $ker(\pi^p_*)$ is acyclic implying
  $\pi^p_*$ is a weak equivalence.  For $p=0$ we have $J^0_*\cong
  K\otimes {\rm CB}_*(C/K)$ and for $p>n$ we have $J^p_n={\rm
  C}_n(K,C/K)$.  The result follows from the fact that $C/K$ is
  H--counital.
\end{proof}

\begin{thm}\label{Hochschild2}
  Let $0\xra{}K\xra{}C\xra{}C/K\xra{}0$ be a short exact sequence of
  H--counital $B$--module coalgebras.  Assume that $X$ has finite
  projective dimension and both $C$ and $C/K$ are projective
  $B$--modules.  Then the canonical inclusion of differential graded
  $k$--modules $k\eotimes{B}{\rm CH}^{\rm tw}_*(K;X)\xra{}
  k\eotimes{B}{\rm CH}^{\rm tw}_*(C,K;X)$ is a weak equivalence.
\end{thm}

\begin{proof}
  Define a sequence $\{F^p_*\}_{p\geq 0}$ of differential graded
  $k$--submodules of ${\rm CH}^{\rm tw}_*(C,K;X)$ by letting
  \begin{align*}
    F^p_n = \begin{cases}
            X\otimes K^{\otimes n+1} & \text{ if } n\leq p\\
            X\otimes K\otimes C^{\otimes n-p}\otimes K^{\otimes p} & \text{ if } n>p
        \end{cases}
  \end{align*}
  There are inclusions $F^{p+1}_*\subseteq F^p_*$ and we see that
  \begin{align*}
    F^p_n/F^{p+1}_n
      = \begin{cases}
    0 & \text{ if } n\leq p\\
        X\otimes K\otimes C^{\otimes n-p-1}\otimes C/K\otimes K^{\otimes p}
          & \text{ if } n>p
    \end{cases}
  \end{align*}
  which is isomorphic to the differential graded $B$--module
  $X\otimes{\rm C}_*(K,C/K)[-p-1]\otimes K^{\otimes p}$.  Since
  $K\xra{}C$ is $B$--split we have a short exact sequence of
  differential graded $k$--modules of the form
  \begin{align*}
    0\xra{}k\eotimes{B}F^{p+1}_*\xra{}k\eotimes{B}F^p_*\xra{}k\eotimes{B}(F^p_*/F^{p+1}_*)\xra{}0
  \end{align*}
  Now use Lemma~\ref{FiniteDimension} to conclude that
  $k\eotimes{B}(F^p_*/F^{p+1}_*)$ is acyclic.  The result follows
  after observing $F^0_*={\rm CH}^{\rm tw}_*(C,K;X)$ and for $n>p$ and
  $F^p_n={\rm CH}^{\rm tw}_n(K;X)$,
\end{proof}

\begin{defn}
  Assume $X$ is a left $B$--module $B\otimes X\xra{\mu_X}X$ and a left
  $B$--comodule $X\xra{\rho_X}B\otimes X$. We will call $X$ stable if
  $\mu_X\rho_X = id_X$.  In other words $x_{(-1)}x_{(0)}=x$ for any
  $x\in X$.
\end{defn}

\begin{exm}
  \begin{enumerate}

    \item Assume $X$ is a $B$--comodule.  One can force $X$ to be a
    left $B$--module via the counit $B\xra{\varepsilon}k$,
    i.e. $b(x)=\varepsilon(b)x$ for any $x\in X$ and $b\in B$.  Then
    $X$ is a stable $B$--module/comodule.  We will denote this stable
    module/comodule by $X_\varepsilon$.

    \item Similarly, if $X$ is a $B$--module, one can force $X$ to be a
    left $B$--comodule via $\rho(x)=(\B{I}\otimes x)$.  Then $X$ is a
    stable $B$--module/comodule.  We will denote this stable
    module/comodule by $X_\B{I}$.

    \item Assume $B$ is a Hopf algebra with an invertible antipode.
    We will consider $B$ as a left $B$--module via the left regular
    representation $b(x)=bx$ and a left $B$--comodule via the left
    co-adjoint coaction $x_{(-1)}\otimes x_{(0)}=
    x_{(1)}S^{-1}(x_{(3)})\otimes x_{(2)}$.  We will denote this
    $B$--module/comodule by $B^{\rm r,ad}$. Then
    \begin{align*}
      x_{(-1)}x_{(0)} = x_{(1)}S^{-1}(x_{(3)})x_{(2)} = x
    \end{align*}
    for any $x\in B$, i.e. $B^{r,ad}$ is a stable $B$--module
    comodule.   Moreover
    \begin{align*}
      (yx)_{(-1)}\otimes (yx)_{(0)}
      = & y_{(1)}x_{(1)}S^{-1}(x_{(3)})S^{-1}(y_{(3)})\otimes y_{(2)}x_{(2)}\\
      = & y_{(1)}x_{(-1)}S^{-1}(y_{(3)})\otimes y_{(2)}x_{(0)}
    \end{align*}
    for any $x,y\in B$ meaning $B^{\rm r,ad}$ is a stable
    anti-Yetter-Drinfeld (SaYD) module.

    \item Assume $B$ is a Hopf algebra with an invertible antipode.
    Consider $B$ as a $B$--module via the left adjoint action $b(x) =
    b_{(1)}xS^{-1}(b_{(2)})$ and as a $B$--comodule via the left
    regular co-representation $x_{(-1)}\otimes x_{(0)}):=
    x_{(1)}\otimes x_{(2)}$.  We will denote this $B$--module comodule
    by $B^{\rm ad,r}$. Then
    \begin{align*}
      x_{(-1)} x_{(0)} = x_{(1)}x_{(3)}S^{-1}(x_{(2)}) = x
    \end{align*}
    for any $x\in B$.  This means $B^{ad,r}$ is a stable
    $B$--module/comodule.  Moreover, $B^{ad,r}$ is also a SaYD module
    since
    \begin{align*}
      (yx)_{(-1)}\otimes (yx)_{(0)}
      = & (y_{(1)}xS^{-1}(y_{(2)}))_{(1)}\otimes (y_{(1)}xS^{-1}(y_{(2)}))_{(2)}\\
      = & y_{(1)}x_{(1)}S^{-1}(y_{(4)})\otimes y_{(2)}x_{(2)}S^{-1}(y_{(3)})\\
      = & y_{(1)}x_{(-1)}S^{-1}(y_{(3)})\otimes y_{(2)}x_{(0)}
    \end{align*}
    for any $x,y\in B$.
  \end{enumerate}
\end{exm}

\begin{rem}
  In \cite{Kaygun:BialgebraCyclicK} we introduced a para-cocyclic
  $B$--module $\B{PCM}^c_*(C,B,X)$ and a cocyclic module
  $\B{CM}^c_*(C,X):= k\eotimes{B} \B{PCM}^c_*(C,B,X)$ of a triple
  $(C,B,X)$ where $B$ is a bialgebra, $C$ is a $B$--module algebra and
  $X$ is a stable $B$--module/comodule.  We also showed that in case
  $B$ is a Hopf algebra with an invertible antipode and $X$ is a
  stable aYD module then one has
  \begin{align*}
    \C{C}^B_*(C;X) \cong \B{CM}^c_*(C,B,X) :=
           k\eotimes{B}\B{PCM}^c_*(C,B,X)
  \end{align*}
  where $\C{C}^B_*(C;X)$ is the cocyclic complex of
  \cite{Khalkhali:HopfCyclicHomology} associated with the coalgebra
  $B$--module coalgebra $C$.  $\C{C}^B_*(C;X)$ is a cocyclic complex
  whose Hochschild complex is $k\eotimes{B}{\rm CH}^{\rm tw}_*(C;X)$
  and whose bar complex is $k\eotimes{B}(X\otimes{\rm CB}_*(C))$ with
  induced differentials.
\end{rem}

\begin{defn}
  A sequence of differential graded (resp.  simplicial, cosimplicial,
  cyclic, cocyclic) modules $\C{X}_*\xra{u_*}\C{Y}_*\xra{v_*}\C{Z}_*$
  is called a homotopy cofibration sequence if there are connecting
  morphisms $H_{n+1}(\C{Z}_*)\xra{\delta_n}H_n(\C{X}_*)$
  (resp. $H_n(\C{Z}_*)\xra{\delta_n}H_{n+1}(\C{X}_*)$ for cosimplicial
  and cocyclic modules) which fit into a long exact sequence of the
  form
  \begin{align*}
    \cdots\xra{H_{n+1}(v_*)}H_{n+1}(\C{Z}_*)\xra{\delta_n}
    H_n(\C{X}_*)\xra{H_n(u_*)}H_n(\C{Y}_*)\xra{H_n(v_*)}
    H_n(\C{Z}_*)\xra{\delta_n}\cdots
  \end{align*}
  or respectively for cosimplicial and cocyclic modules of the form
  \begin{align*}
    \cdots\xra{H_{n-1}(v_*)}H_{n-1}(\C{Z}_*)\xra{\delta_{n-1}}
    H_n(\C{X}_*)\xra{H_n(u_*)}H_n(\C{Y}_*)\xra{H_n(v_*)}
    H_n(\C{Z}_*)\xra{\delta_n}\cdots
  \end{align*}
  for any $n\geq 0$.
\end{defn}

\begin{thm}[Excision in Hopf Cyclic Homology]\label{Excision1}
  Let $B$ be a Hopf algebra with an invertible antipode and $X$ be a
  stable anti-Yetter-Drinfeld (SaYD) module.  Let $C$ be a counital
  $B$--module coalgebra and $0\xra{}K\xra{}C\xra{}C/K\xra{}0$ be a
  short exact sequence of H--counital $B$--module coalgebras.  Assume
  that both $C$ and $C/K$ are $B$--projective and $X$ has finite
  projective dimension.  Then there is a homotopy cofibration sequence
  of the form
  \begin{align*}
    \B{CM}^c_*(K,B,X)\xra{}\B{CM}^c_*(C,B,X)\xra{}\B{CM}^c_*(C/K,B,X)
  \end{align*}
\end{thm}

\begin{proof}
  Theorem~\ref{Hochschild0}, Theorem~\ref{Hochschild1} and
  Theorem~\ref{Hochschild2} tell us we have a homotopy cofibration
  sequence of the form
  \begin{align}\label{ExcisionDiagram}
    k\eotimes{B}{\rm CH}^{\rm tw}_*(K;X)\xra{}
    k\eotimes{B}{\rm CH}^{\rm tw}_*(C;X)\xra{}
    k\eotimes{B}{\rm CH}^{\rm tw}_*(C/K;X)
  \end{align}
  This is equivalent to saying Hopf Hochschild homology has excision.
  Since $K$, $C$ and $C/K$ are all H--counital and $B$--projective,
  the differential graded $k$--complexes $k\eotimes{B}(X\otimes{\rm
  CB}_*(K))$, $k\eotimes{B}(X\otimes{\rm CB}_*(C))$ and
  $k\eotimes{B}(X\otimes{\rm CB}_*(C/K))$ are all acyclic because $X$
  has finite projective dimension following
  Lemma~\ref{FiniteDimension}.  Therefore we have another homotopy
  cofibration sequence
  \begin{align*}
    k\eotimes{B}(X\otimes{\rm CB}_*(K))\xra{}
    k\eotimes{B}(X\otimes{\rm CB}_*(C))\xra{}
    k\eotimes{B}(X\otimes{\rm CB}_*(C/K))
  \end{align*}
  Since $\B{CM}^c_*(\cdot;X)$ is composed of $k\eotimes{B}{\rm
  CH}^{\rm tw}_*(\cdot;X)$ and $k\eotimes{B}(X\otimes{\rm
  CB}_*(\cdot))$ with induced cocyclic structures the result follows.
\end{proof}

\begin{cor}
  Assume $B$ is a Hopf algebra with an invertible antipode and $X$ be
  a SaYD module.  Let $\{C_i\}_{i\in \Lambda}$ be a class of counital
  $B$--module coalgebras such that each $C_i$ is $B$--projective.  Let
  $C=\bigoplus_{i\in\Lambda}C_i$ be the direct sum considered as a
  counital $B$--module coalgebra via the direct sum of
  comultiplication structures.  Then
  \begin{align*}
    HC_*(C;X)\cong \bigoplus_{i\in\Lambda} HC_*(C_i;X)
  \end{align*}
  In other words cyclic homology is additive on the class of
  $B$--projective counital $B$--module coalgebras.
\end{cor}

\begin{rem}
  If we drop $B$, in other words if we set $B=k$, then
  Theorem~\ref{Excision1} reduces to a Wodzicki type excision theorem
  \cite{FarinatiSolotar:ExcisionInCyclicHomologyOfCoalgebras} for
  coalgebras.  We suspect that we have a Cuntz--Quillen type of
  excision \cite{QuillenCuntz:Excision} in periodic equivariant cyclic
  homology of coalgebras, if we drop the H-counitality condition.
\end{rem}

\begin{conj}\label{Conjecture1}
  Let $B$ be a Hopf algebra with an invertible antipode and $X$ be a
  SaYD module.  Let $0\xra{}K\xra{}C\xra{}C/K\xra{}0$ be a sequence of
  H--counital $B$--module coalgebras such that both $C$ and $C/K$ are
  $B$--projective.  Then one has a 6-term exact sequence in periodic
  cyclic homology of the form
  \begin{align*}
    \vspace{5mm}\\
    \begin{CD}
    HP_0(K,B;X)   @>>> HP_0(C,B;X) @>>> HP_0(C/K,B;X)\\
    @AAA             @.                 @VVV\\
    HP_1(C/K,B;X) @<<< HP_1(C,B;X) @<<< HP_1(K,B;X)
    \end{CD}\\
  \end{align*}
\end{conj}

\section{Co-integrals and finite projective dimension}

The following proposition gives us a good example of a Hopf algebra
and a coefficient module with finite projective dimension.

\begin{prop}
  If $B$ is $U(\G{g})$ the universal enveloping algebra of a finite
  dimensional Lie algebra, then the trivial $\G{g}$--module $k$ has
  finite projective dimension.
\end{prop}

\begin{proof}
  The Chevalley--Eilenberg complex ${\rm
  CE}_*(\G{g})=U(\G{g})\otimes\Lambda^*\G{g}$ of $\G{g}$ with
  coefficients in $U(\G{g})$ is a $U(\G{g})$--projective resolution of
  $k$.  Since $\G{g}$ is finite dimensional, ${\rm CE}_*(\G{g})$ is a
  finite resolution.
\end{proof}

In the rest of this section, we investigate the conditions which would
make $k_{(1,\epsilon)}$ into a projective $B$--module.

Recall that an element $\sigma\in B$ is called a normalized right
co-integral if one has $\sigma b = \varepsilon(b)\sigma$ and
$\varepsilon(\sigma) = 1$ for any $b\in B$.  It is well-known
\cite{Sweedler:HopfAlgebras} that if $B$ is a Hopf algebra then $B$
admits a right (or equivalently a left) co-integral iff $B$ is
semi-simple.

The following fact is easy to prove, yet its proof is given to
demonstrate the connection between the $B$--projectivity of $k$ and
existence of a co-integral.

\begin{prop}
  $k$ is a projective right $B$--module iff $B$ admits a normalized
  right co-integral.
\end{prop}

\begin{proof}
  Assume $k$ is projective right $B$--module.  Then since any
  epimorphism of right $B$--modules of the form $Z\xra{p}k$ splits, so
  does $B\xra{\varepsilon}k$.  Let us denote the splitting by
  $k\xra{s}B$.  Therefore there is an element $\sigma=s(1)$ such that
  $\sigma b=s(1)b=s(\varepsilon(b) 1)=e=\varepsilon(b)\sigma$ for any $b\in
  B$ implying $\sigma$ is a right co-integral Conversely, assume $B$
  admits a right co-integral $\sigma\in B$ and consider an arbitrary
  epimorphism of right $B$--modules of the form $Z\xra{p}k$.  Pick an
  element $z\in Z$ such that $p(z_1)=1$.  Since $p$ is non-trivial and
  $k$ is a field, such an element exists.  Define $k\xra{s}Z$ by
  letting $s(\lambda)=\lambda z_1\sigma$ for any $\lambda\in k$.  Then
  $s$ is a morphism of $B$--modules since
  \begin{align*}
    s(\lambda b)=\lambda\varepsilon(b)s(1)=\lambda\varepsilon(b)z_1\sigma = \lambda
     z_1\sigma b = s(\lambda) b
  \end{align*}
  Moreover,
  \begin{align*}
    p(s(\lambda)) = \lambda p(z_1\sigma) = \lambda p(z_1)\varepsilon(\sigma)
    = \lambda
  \end{align*}
  implying $s$ is a right $B$--module splitting of $p$.  Result
  follows.
\end{proof}

\begin{cor}
  If $B$ is a Hopf algebra admitting a normalized right or left
  co-integral and $C$ is any counital $B$--module coalgebra, then for
  any SaYD module $X$ and for any short exact sequence of H--counital
  $B$--module coalgebras $0\xra{}K\xra{}C\xra{}C/K\xra{}0$ one has a
  homotopy cofibration sequence given in
  Equation~\ref{ExcisionDiagram}
\end{cor}

\section{Relative bialgebra cyclic homology}

\begin{defn}
  Given a $B$--module coalgebra $C$ and a $B$--submodule subcoalgebra
  $K$, we define the relative cyclic homology of the pair $(C,K)$ with
  coefficients in a stable $B$--module/comodule $X$ as the cyclic
  homology of the cyclic $k$--module $\B{CM}^c_*((C,K),B,X)$ which is
  defined as
  \begin{align*}
    coker\big(\B{CM}^c_*(K,B,X)\xra{\B{CM}^c_*(i)}\B{CM}^c_*(C,B,X)\big)
  \end{align*}
  the cokernel of the canonical morphism of cocyclic $k$--modules
  $\B{CM}^c_*(K,B,X)\xra{\B{CM}^c_*(i)}\B{CM}^c_*(C,B,X)$ coming from
  the inclusion $K\xra{}C$.
\end{defn}

\begin{rem}
  The canonical morphism $\B{CM}^c_*(K,B,X)\xra{}\B{CM}^c_*(C,B,X)$ is
  not necessarily a monomorphism of cocyclic $k$--modules.  When
  $K\xra{}C$ is $B$--split as a morphism of $B$--modules this
  canonical inclusion is a monomorphism of cocyclic $k$--modules.
\end{rem}

\begin{defn}
  Given a $B$--module coalgebra $C$ and a $B$--submodule coideal $J$,
  we define the relative cyclic homology of the pair $(C,J)$ with
  coefficients in a stable $B$--module/comodule $X$ as the cyclic
  homology of the triple $(C/J,B,X)$.
\end{defn}

\begin{rem}
  Assume $C$ is a $B$--module coalgebra.  Any $B$--submodule
  subcoalgebra $K$ is also a $B$--submodule coideal which means
  there are potentially two different versions of relative cyclic
  homology of the pair $(C,K)$.
\end{rem}

\begin{rem}
  In \cite{ConnesMoscovici:HopfRelativeCyclic} Connes and Moscovici
  defined the relative Hopf cyclic homology of a Hopf algebra $H$ and
  a Hopf subalgebra $K$ with coefficients in a SaYD module $X$ as the
  Hopf cyclic homology of the triple $(H\eotimes{K}k,H,X)$.  This is
  the same as the relative cyclic homology of the pair $(H,HK^+)$ with
  coefficients in $X$ where $HK^+$ is the left ideal of $H$ generated
  by $ker(\varepsilon_K)$ which is also a $H$--submodule coideal of
  $H$.
\end{rem}

\begin{thm}
  Assume $B$ is a Hopf algebra with an invertible antipode and $X$ is
  a SaYD module.  Assume also that $X$ has finite projective
  $B$--dimension.  Now, let $0\xra{}K\xra{}C\xra{}C/K\xra{}0$ be a
  short exact sequence of H--counital $B$--module coalgebras such that
  both $C$ and $C/K$ are $B$--projective.  Then two relative cyclic
  homologies we defined above for the pair $(C,K)$ with coefficients
  in $X$ are isomorphic.
\end{thm}

\begin{proof}
  Our definition dictates that the relative cyclic homology is defined
  through the short exact sequence 
  \begin{align*}
    0\xra{}\B{CM}^c_*(K,C,X)\xra{}\B{CM}^c_*(C,B,X)\xra{}\B{CM}^c_*((C,K),B,X)\xra{}0
  \end{align*}
  However we also have the homotopy cofibration sequence
  \begin{align*}
    \B{CM}^c_*(K,B,X)\xra{}\B{CM}^c_*(C,B,X)\xra{}\B{CM}^c_*(C/K,B,X)
  \end{align*}
  meaning the two different relative theories we defined above are
  isomorphic.
\end{proof}

\begin{defn}
  A $k$--submodule $J$ of a bialgebra $B$ is called a bialgebra ideal
  if (i) $J$ is a two sided ideal of $B$ (ii) $J$ is a two sided
  coideal of $B$ (iii) $\varepsilon(J)\equiv 0$.  If $B$ is a Hopf
  algebra then $J$ is called a Hopf ideal if $J$ is a bialgebra ideal
  and $S(J)=J$.
\end{defn}

\begin{thm}\label{CommutativeHopf}
  Let $B$ be a bialgebra and $J$ be a bialgebra ideal of $B$.  Assume
  also that $X$ is a stable $B$--module/comodule such that $z\cdot
  x=0$ for any $z\in J$.  Then the relative cyclic homology of the
  pair $(B,J)$ with coefficients in $X$ is the same as the cyclic
  homology of the triple $(B/J,B/J,X)$.
\end{thm}

\begin{proof}
  Since $X$ is annihilated by $J$, the action of $B$ on $X$ actually
  splits as
  \begin{align*}
    B\otimes X\xra{q_{B/J}\otimes id_X}B/J\otimes X\xra{\mu_X}X
  \end{align*}
  where $B\xra{q_{B/J}}B/J$ is the canonical quotient morphism.  Thus $X$ is a
  stable $B/J$--module/comodule.  Moreover,
  \begin{align*}
    k\eotimes{B}\left((B/J)^{\otimes n+1}\otimes X\right)
    = k\eotimes{B/J}\left((B/J)^{\otimes n+1}\otimes X\right)
  \end{align*}
  for any $n\geq 0$ which implies
  \begin{align*}
    \B{CM}^c_*(B/J,B,X) = \B{CM}^c_*(B/J,B/J,X)
  \end{align*}
  as we wanted to show.
\end{proof}

\begin{cor}
  Let $B$ be a commutative Hopf algebra and $J$ be a Hopf ideal of
  $B$.  Then the relative cyclic homology of the pair $(B,J)$ with
  coefficients in the 1-dimensional stable $B$--module/comodule $k$ is
  \begin{align*}
    {\rm HC}^{\B{CM},c}_n((B,J),B,k)\cong 
    \bigoplus_{i\geq 0}{\rm HH}_{n-2i}(B/J,k)
  \end{align*}
  for any $n\geq 0$.
\end{cor}

\begin{proof}
  Theorem~\ref{CommutativeHopf} implies $\B{CM}^c_*(B/J,B,X)\cong
  \B{CM}^c_*(B/J,B/J,X)$.  Then Theorem 4.2 of
  \cite{Khalkhali:DualCyclicHomology} implies the result.
\end{proof}

\section{The dual theory}

In this section we will abstain from giving detailed proofs since they
can be obtained from carefully dualizing corresponding proofs for the
Hopf cyclic homology of module coalgebras given above.

\subsection{Definitions}

  Recall that an an algebra $A$ is called a right $B$--comodule
  algebra if $A$ is a right $B$--comodule $A\xra{\rho_A}A\otimes B$
  and for any $a^1,a^2\in A$ we have
  \begin{align*}
    (a^1a^2)_{(0)}\otimes (a^1a^2)_{(1)}
    = a^1_{(0)}a^2_{(0)}\otimes a^1_{(1)}a^2_{(1)}
  \end{align*}
  which is equivalent to saying $\rho_A$ is a morphism of algebras.

  In \cite{Kaygun:BialgebraCyclicK} we constructed a cyclic theory
  ${\rm HC}^{\B{CM},a}_*(A,B,X)$ where $B$ is bialgebra, $A$ is a
  $B$--comodule algebra and $X$ is just a stable $B$--module/comodule.
  We also showed that in case $B$ is a Hopf algebra with an invertible
  antipode and $X$ is a stable anti-Yetter-Drinfeld module the theory
  we developed and the theory developed in
  \cite{Khalkhali:HopfCyclicHomology} are the same. Here we develop a
  relative version of this dual theory in the bialgebra setting.

\begin{defn}
  Given a $B$--comodule algebra $A$ and a $B$--subcomodule bilateral
  ideal $I\subset A$, we define the dual relative Hopf cyclic homology
  of the pair $(A,I)$ with coefficients in a stable
  $B$--module/comodule $X$ as the dual cyclic homology of the triple
  $(A/I,B,X)$.
\end{defn}

\begin{defn}
  Given a $B$--comodule algebra $A$ and a $B$--subcomodule subalgebra
  $A' \subset A$ we define the relative Hopf cyclic homology of
  $(A,A')$ with coefficients in a stable $B$--module/comodule $X$ as
  the homology of the complex
  \begin{align*}
    coker\big(\B{CM}^a_*(A',B,X)\xra{\B{CM}^a_*(i)}\B{CM}^a_*(A,B,X)\big)
  \end{align*}
  the cokernel of the canonical morphism of cyclic objects
  $\B{CM}^a_*(A',B,X)\xra{\B{CM}^a_*(i)}\B{CM}^a_*(A,B,X)$ coming from
  the inclusion $A'\xra{i}A$. This cyclic object is denoted by
  $\B{CM}^a_*((A,A'),B,X)$.
\end{defn}

\begin{rem}
  Since any bilateral ideal is also a subalgebra, there are two
  possibly different definitions of dual relative cyclic homology.
  These different definitions will agree when the canonical morphism
  of cyclic objects $\B{CM}^a_*((A,I),B,X)\xra{}\B{CM}^a_*(A/I,B,X)$
  is a weak equivalence.  This will happen if and only if we have
  excision in dual bialgebra cyclic homology.
\end{rem}

\subsection{Cocommutative Hopf algebras}

\begin{thm}\label{CocommutativeHopf}
  Assume $K$ is a bialgebra ideal of $B$ such that the canonical
  epimorphism of coalgebras $B\xra{q} B/K$ splits in the category of
  coalgebras.  Assume also that $X$ is a stable $B$--module/comodule
  with the property that ${\rm Res}^B_K(X)= K\fotimes{B}X = 0$.  Then
  the dual relative cyclic homology of the pair $(B,K)$ with
  coefficients in $X$ is the same as the dual cyclic homology of the
  triple $(B/K,B/K,X)$.
\end{thm}

\begin{proof}
  Since ${\rm Res}^B_K(X)=0$ and $B\xra{q}B/K$ splits, the $B$--coaction
  on $X$ splits as
  \begin{align*}
    X\xra{\rho_X}X\otimes B/K\xra{}X\otimes B
  \end{align*}
  This means we have an equality
  \begin{align*}
    \left((B/K)^{\otimes n+1}\otimes X\right)\fotimes{B}k
    = \left((B/K)^{\otimes n+1}\otimes X\right)\fotimes{B/K}k
  \end{align*}
  which means $\B{CM}^a_*(B/K,B,X)=\B{CM}^a_*(B/K,B/K,X)$.  The result 
follows.
\end{proof}

\begin{cor}
  Let $B$ be a cocommutative Hopf algebra and $K$ be a Hopf ideal of
  $B$ such that the epimorphism $B\xra{q}B/K$ splits as a morphism of
  coalgebras.  Then one has
  \begin{align*}
    {\rm HC}^{\B{CM},a}_n(B/K,B,k_{(1,\varepsilon)}) \cong
    \bigoplus_{i\geq 0} {\rm HH}_{n-2i}(B/K,k)
  \end{align*}
\end{cor}

\begin{proof}
  Theorem~\ref{CocommutativeHopf} above and Theorem~4.1 of
  \cite{Khalkhali:DualCyclicHomology}.
\end{proof}

\begin{exm}
  Let $G$ be an arbitrary discrete group and $H$ be a normal subgroup
  of $G$.  Consider the two sided ideal $a_H$ of $k[G]$ generated by
  elements of the form $(h-1)$ where $h\in H$.  First, observe that
  $a_H$ is also a two sided coideal of the Hopf algebra $k[G]$ and is
  stable under the antipode map $S(g)=g^{-1}$.  Next observe that
  $k[G]/a_H\cong k[G]\eotimes{H}k\cong k[G/H]$ and
  $a_H\fotimes{G}k=0$.  Now for every coset $gH$ in $G/H$ pick a
  representative $s(gH)\in G$.  The induced map $k[G/H]\xra{s}k[G]$ is
  not a morphism of algebras but is a morphism of coalgebras since
  $\Delta(x)=(x\otimes x)$ for any $x$ from $k[G]$.  So, the canonical
  quotient $k[G]\xra{}k[G/H]$ splits as a morphism of coalgebras.
  Therefore the relative Hopf cyclic homology of $k[G]$ with respect
  to the Hopf ideal $a_H$ with coefficients in the trivial
  $k[G]$--module/comodule $k$ is given by
  \begin{align}
    HC^{\B{CM},a}_n(k[G]/a_H,k[G];k)
    \cong \bigoplus_{i\geq 0}HH_{n-2i}(k[G/H],k)
    \cong \bigoplus_{i\geq 0}H_{n-2i}(G/H,k)
  \end{align}
\end{exm}

\subsection{Excision in dual Hopf cyclic homology}

\begin{thm}[Excision in Dual Hopf Cyclic Homology]\label{Excision2}
  Let $B$ be a Hopf algebra with an invertible antipode and $X$ be a
  stable anti-Yetter-Drinfeld module.  Assume $A$ is a unital
  $B$--comodule algebra and $I$ is a H--unital $B$-subcomodule ideal
  of $A$.  Assume also that both $I$ and $A$ are
  $B$--co-projective. Then the short exact sequence of $B$--comodule
  algebras $0\xra{}I\xra{}A\xra{}A/I\xra{}0$ yields a homotopy
  cofibration sequence
  \begin{align}\label{HomotopyCofibration}
    \B{CM}^a_*(I,B,X)\xra{}\B{CM}^a_*(A,B,X)\xra{}\B{CM}^a_*(A/I,B,X)
  \end{align}
\end{thm}

\begin{rem}
  In this version of the excision, we do not need any condition on the
  $B$--co-projective dimension of $X$ since the resulting spectral
  sequence computing the hyper Cotor functors is convergent even when
  $X$ does not have finite co-projective dimension.
\end{rem}

\begin{rem}
  If we drop $B$, in other words if we set $B=k$, then
  Theorem~\ref{Excision2} reduces to Wodzicki's excision theorem
  \cite{Wodzicki:Excision}.  Again, we suspect that we have a
  Cuntz--Quillen type excision \cite{QuillenCuntz:Excision} in
  periodic equivariant cyclic homology, if we drop the H-unitality
  condition.
\end{rem}

\begin{conj}\label{Conjecture2}
  Let $B$ be a Hopf algebra with an invertible antipode and $X$ be a
  SaYD module.  Let $0\xra{}J\xra{}A\xra{}A/J\xra{}0$ be a sequence of
  $B$--comodule algebras such that both $J$ and $A$ are
  $B$--co-projective.  Then one has a 6-term exact sequence in
  periodic cyclic homology of the form
  \begin{align*}
    \vspace{3mm}\\
    \begin{CD}
    HP_1(J,B;X)   @>>> HP_1(A,B;X) @>>> HP_1(A/J,B;X)\\
    @AAA             @.                 @VVV\\
    HP_0(A/J,B;X) @<<< HP_0(A,B;X) @<<< HP_0(J,B;X)
    \end{CD}\\
  \end{align*}
\end{conj}

\subsection{Integrals}

  Recall that a linear functional $B\xra{\eta}k$ is called a
  normalized left integral if
  \begin{align*}
    b_{(1)}\eta(b_{(2)}) = \eta(b)\B{I}
    \text{\ \ \ \ and\ \ \ \ }
    \eta(\B{I})  = 1
  \end{align*}
  for any $b\in B$.  Right integrals are defined similarly.

  Let $B$ be a Hopf algebra with an invertible antipode.  Then it is
  easy to see that $B$ admits a normalized left integral iff $B$
  admits a normalized right integral.  The following theorem is
  standard (cf. e.g. \cite{Sweedler:HopfAlgebras} for a proof):

\begin{prop}
  Assume $B$ is a Hopf algebra with an invertible antipode.  $B$
  admits a left integral iff $B$ is co-semi-simple as a coalgebra.
\end{prop}

\begin{cor}
  Let $B$ be a Hopf algebra with an invertible antipode and $X$ be a
  SaYD module.  Assume $B$ admits a left integral.  Let $A$ be a
  unital $B$--comodule algebra and $I$ be a H--unital $B$--submodule
  ideal of $A$.  Then there is a homotopy cofibration sequence given
  in Equation~\ref{HomotopyCofibration}.
\end{cor}

\begin{proof}
  Existence of left integral implies $B$ is co-semi-simple which in
  turn implies every $B$--comodule is $B$--co-projective.  The result
  follows.
\end{proof}


\end{document}